\documentclass{article}    
\usepackage{latexsym}

 \newcommand{\N}{{\rm I~\hspace{-1.10ex}N}} 
\newcommand{\R}{{\rm I~\hspace{-1.15ex}R}}

\newtheorem{Theorem}{Theorem}[section]
\newtheorem{Lemma}{Lemma}[section]

\newtheorem{remark}{Remark}[section]

\newtheorem{definition}{Definition}[section]

\title{On a  stochastic partial differential  equation with non-local diffusion.} 
\author{Pascal AZERAD
and Mohamed MELLOUK,\footnote{\texttt{\{azerad,mellouk\}@math.univ-montp2.fr}}\medskip\\
Institut de Math\'ematiques et de Mod\'elisation de Montpellier,\\
UMR 5149, Universit\'e Montpellier 2, cc 51\\
34095  Montpellier Cedex 5, France.
}
\date{\today}
\begin{document}                                                                \maketitle

\begin{abstract}
In this paper, we prove existence, uniqueness and regularity 
 for a class
of  stochastic partial
differential equations  with a fractional Laplacian driven by a space-time white noise in dimension one. The equation we consider may also include a reaction term.
\end{abstract}

\section{Introduction and general framework}
In recent years, fractional calculus has received a great deal of attention. Equations involving fractional
derivatives and fractional Laplacians have been studied by various authors (see, e.g. Podlubny 
 \cite{Po} and references
therein). In probability theory, fractional calculus has been extensively used in the study of fractional 
Brownian motions. In this work we  consider a stochastic partial differential equations where the  standard Laplacian
operator is replaced by a fractional one. \\

\noindent
Let $\lambda > 0$. We consider the fractional Laplacian 
$\Delta_\lambda=- (-\Delta)^{\lambda/2}= - (-\partial^2/\partial x^2)^{\lambda/2}$, the symmetric fractional derivative 
 of order $\lambda$ on $\R$. This is a non-local operator defined  via the Fourier transform
 ${\cal F}$:
$$ {\cal F}(\Delta_\lambda v)(\xi)=-\vert \xi\vert^\lambda {\cal F}(v)(\xi).$$
It also has another representation, for $0<\lambda <2$,
\begin{equation}
\label{delta}
\Delta_\lambda v(x)=K \int_{\R} \left\{ v(x+y)-v(x)-
\nabla v(x)\cdot \frac{y}{1+\vert y\vert^2}
\right\}\frac{dy}{\vert y\vert^{1+\lambda}} ,
\end{equation}
 for some positive constant $K=K_\lambda$, which identifies it as the infinitesimal generator for
 the symmetric 
 $\lambda$-stable L\'evy process (see, e.g., It\^o \cite{It}, Stroock \cite{St}, Komatsu \cite{Ko},
 Dawson and Gorostiza \cite{Da}).\\

\noindent
 
 Let $W=\left\{W(t,x),(t,x)\in [0,T]\times \R \right\}$ be a Brownian sheet on a complete probability
 space $(\Omega ,{\cal{G}},P)$. That is, $W$ is a zero-mean Gaussian random field with covariance
 function
$$
E(W(t,x)W(s,y))=\frac{1}{2} (
s\wedge t) \left( \vert x\vert +\vert y\vert -\vert x-y\vert\right), 
$$
$x, y \in \R$, $s, t \in [0,T]$. 
Then, for each $t\in [0,T]$, we define a filtration
$$
{\cal{G}}_t^0=\sigma \left( W(s,x), s\in [0,t], x\in \R\right),\,\,\,\,\,
 {\cal{G}}_t= {\cal{G}}_t^0 \vee {\cal{N}}, 
$$
where ${\cal{N}}$ is the $\sigma$-field generated by sets with  $P$\--outer measure zero. \\

The family of $\sigma$-fields  $\{{\cal{G}}_t, 0\le t\le T\}$ constitutes a stochastic basis on the probability
space  $(\Omega ,{\cal{G}},P)$. Let ${\cal P}$ the corresponding predictable $\sigma$-field on
$\Omega \times [0,T] \times \R$. The stochastic integral with respect to the Brownian sheet
is explained in Cairoli et  al. \cite{Ca} or Walsh \cite{Wa}.\\

We focus on the  following  parabolic stochastic partial
differential equation, driven by space--time white noise in one space
dimension on $[0, T]\times \R $%
$${\bf (E)}\hspace{4mm}\frac{\partial u}{\partial t}(t,x)=\Delta_\lambda u(t,x)
+b\left( t,x,u\left( t,x\right) \right)+\sigma \left( t,x,u\left( t,x\right) \right)\dot W(t,x),$$
with initial condition $u(0,x)=u_0(x)$  ${\cal{G}}_0$-measurable and satisfying some conditions that will be specified
later. The process $\dot W(t,x) = \frac{\partial^2 W}{\partial t \partial x}$ is 
the generalized (distribution) derivative of the Brownian sheet. 
The properties of $\dot W$ are described in Walsh \cite{Wa}.\\

In principle one can think of a wide variety of random forcing terms. White noise in time and space
is very often a candidate. Main motivations behind this choice are  
central limit type theorems and the insufficient knowledge of the
 neglected effects or external disturbances.\\

Evolution problems involving fractional Laplace operator have long been extensively studied in mathematical and physical
literature. In the latter, this type of models has been motivated by fractal (anomalous) diffusion related to the L\'evy flights 
(see, e.g., Stroock \cite{St}, Bardos et al. \cite{Ba}, Dawson and Gorostiza
\cite{Da}, Metzler and Klafter \cite{Me}, Mann and
Woyczynski \cite{Man}). In fact, in various physical 
phenomena in statistical mechanics, 
 the anomalous diffusive terms can be nonlocal and fractal, i.e.  
represented by a fractional power of the Laplacian.\\

Equation {\bf (E)} is a generalization of the classical stochastic heat equation where $\lambda=2$ (see, e.g.,
Walsh \cite{Wa}, Pardoux \cite{Pa2}  and the references quoted therein).
 In those papers, the authors
prove  existence and uniqueness of the mild solution in the space interval $[0,1]$. 
The proof relies stronly  on  properties of the explicit  Green kernel associated 
to the operator $\frac{\partial^2}{\partial x^2}$ in   bounded space interval 
with Dirichlet boundary conditions. In  the present paper, we consider the above
 class of equations in the whole line, instead  of a bounded interval, 
for the space variable. The main  properties of  the semigroup generated by
 the fractional Laplacian can be derived by Fourier transform techniques.\\

Consider the fundamental solution  $G_\lambda(t,x)$, associated to the equation
{\bf (E)}  on $[0,T] \times \R $ i.e.
 the convolution kernel of the L\'evy semigroup $\exp(t \Delta_\lambda)$ in $\R$.\\

Using  Fourier transform, we easily see that  $G_\lambda(t,x)$ is given by :
$$G_\lambda(t,x)= {{\cal F}}^{-1}(e^{-t \vert\,\cdot\,\vert^\lambda})(x)
={\int_{\R}} e^{2 i\pi  x\xi} e^{-t\vert\xi\vert^\lambda}d\xi={\cal{F}}(e^{-t 
\vert\,\cdot\,\vert^\lambda})(x).$$
 For $\lambda \in ]0,2]$, the most important property of $G_\lambda$ is its \emph{nonnegativity} (see  L\'evy \cite{Le}
or Droniou et al. \cite{Dro1} for a quick proof ). \\

Throughout this work we consider solutions to the spde {\bf (E)} in the mild sense, 
following Walsh \cite{Wa}, given by the following definition
(which is formally equivalent to  Duhamel's principle or the variation of parameters formula): \\

\begin{definition}
A stochastic process $u: \Omega \times [0,T] \times \R \rightarrow \R$, which is jointly measurable
and ${\cal{G}}_t$-adapted, is said to be a (stochastically) mild solution to the stochastic 
equation {\bf (E)}  with initial condition $u_0$ if there exists a martingale measure $W$, defined on 
$\Omega$, such that 
 {\it a.s.} for almost all $t\in [0,T], x\in \R$,
\begin{eqnarray}
\label{walsh}
u(t,x)&=& G_\lambda(t,\cdot) \ast u_0(x)+\int_0^t\hspace{-2mm}\int_\R
G_\lambda(t-s,x-y)b(s,y,u(s,y))dyds\nonumber\\
 & & \nonumber\\ 
& & \quad  \quad +\int_0^t
\hspace{-2mm}\int_\R G_\lambda(t-s,x-y)\sigma (s,y,u(s,y))W(dy,ds) ,
\end{eqnarray}
 where the last integral is an It\^o stochastic integral.
\end{definition}

\noindent
We assume that the reaction term  $b$ and the white-noise amplitude $\sigma$ are continuous functions on
$[0,T] \times \R \times \R$
  and satisfy the following growth and Lipschitz conditions:\\

\noindent
${\bf (H_0)}$ \\
For all $T>0$, there exists a constants $C=C(T)$, such that for 
all $0\le t \le T, x\in \R$
and $u \in \R$,
\begin{eqnarray*}
 \vert b(t,x,u)\vert + \vert \sigma (t,x,u)\vert & \le & C (1+ \vert u\vert),\\
\vert \sigma (t,x,u)- \sigma (t,x,v)\vert & \le &
 C \,\vert\, u - v\,\vert.\\
\vert b (s,x,u)- b(t,y,v)\vert  & \le &
 C \left( \vert\, t - s\,\vert+\vert\, x - y\,\vert+\vert\, u - v\,\vert \right).
\end{eqnarray*}

\medskip
We shall also need some hypotheses on the initial condition $u_0$ :

\noindent
${\bf (H_1.1)}\quad \sup_{x\in \R} E(\displaystyle|u_0(x)|^p)<\infty,\,\forall p\in[1,+\infty[$.\\

\medskip  
\noindent  
${\bf (H_1.2)}\quad  \exists\rho \in (0,1),\,\forall z \in \R,\,\forall p\in[1,+\infty[,\, \exists C_p>0 $
$$ 
\sup_{y\in \R}E|u_0(y+z)- u_0(y)|^p \leq C_p \vert z \vert^{\rho p}.  
$$  
\bigskip\\
Let us  recall some  well-known  properties 
(see, e.g. Komatsu \cite{Ko}, Biler et Woyczynski \cite{Bi}, 
Droniou et Imbert \cite{Dro})
 of the Green kernel $G_\lambda(t,x)$ which will be used later on.
\begin{Lemma}
\label{tech} Let $\lambda \in ]0,2]$. 
The convolution kernel $G_\lambda$ satisfies the following properties:
\mbox{ }\\
(\textbf{a}) For any $t\in \left] 0,+\infty \right[ $ and $x\in \R $,
$$
G_\lambda(t,x) \geq 0 \quad\mathrm{and}\quad \displaystyle\int_{\R}G_\lambda(t,x)dx = 1. 
$$
(\textbf{b})\ (self similarity) For any $t\in \R_{+}$\ and\ $x\in \R $
$$
G_\lambda(t,x)=t^{-\frac{1}{\lambda}} G_\lambda(1,t^{-\frac{1}{\lambda}} x), 
$$
(\textbf{c}) $G_\lambda$ is $C^\infty$ on $\left] 0,+\infty \right[ \times \R$ and, for 
$m\geq 0$, there exists \\ $C_m>0$ such that for any $t\in \R_{+}$ and $x\in \R $ 
$$
\mid \partial_x^m G_\lambda(t,x)\mid \leq \frac{1}{t^{(1+m)/\lambda}} \frac{C_m}{(1+t^{-2/\lambda}\vert x\vert^2)}. 
$$
(\textbf{d}) For any $(s,t) \in \,]0,\infty[\times ]0,\infty[$
 $$G_\lambda(s,\cdot)\ast G_\lambda(t,\cdot)=G_\lambda(s+t,\cdot).$$ 
(\textbf{e}) $\int_0^T\,dt\int_\R\,dx \, G_\lambda(t,x)^\alpha <\infty$ iff $1/2 < \alpha < 1+\lambda$.
\end{Lemma} 

\noindent
In this paper, in order to  define the stochastic integral, we restrict ourselves to the case
 $\lambda \in ]1,2]$ :  we must take $\lambda \leq 2$ to have $G_\lambda$
positive and  we have to take $\lambda >1$ in order that
 $\int_0^T\int_\R G_\lambda(t,x)^2\,dtdx <\infty$, by  lemma  \ref{tech} (\textbf{e}).\\ 

\noindent
Inessential constants will be denoted
 generically by $C$, even if they vary from line to line.\\

\noindent
The paper is organized as follows. In section 2, we prove 
existence and uniqueness of the solution.
In section 3 we prove  H\"older continuity of the solution in space and time. 
A Gronwall-type improved inequality  and 
an  H\"older inequality frequently used in the paper are 
collected in the appendix.

\section{ Existence and Uniqueness of the solution}
\setcounter{equation}{0}

The main result of this section is the following:

\medskip
\begin{Theorem}
\label{prop1} Let  $\lambda \in ]1,2]$.
Suppose that the hypothesis ${\bf (H_0)}$  and  ${\bf (H_1.1)}$ hold. Then
 there exists a unique solution  $u(t,x)$  to {\bf (E)}  such that: for any $T>0$ and $p\geq 1$,

\begin{equation}
\label{borne0}
\sup_{0\le t\le T}\sup_{x \in \R}E(|u(t,x)|^p)\leq
 C_p<\infty.
\end{equation} 
\end{Theorem}

\medskip



\noindent
{\bf Proof.} The proof of the existence can be done by the usual Picard
iteration procedure. That is, we define recursively 
 
$$u^0\left( t,x\right)  =  \displaystyle \int_{\R} G_\lambda(t, x-y) u_0\left( y\right) dy,$$  
 
\begin{equation}  \label{picard}
\begin{array}{lll}
u^{n+1}\left( t,x\right) & = & u^0\left( t,x\right) + 
\displaystyle\int_0^t\hspace{-2mm}\displaystyle\int_{\R}G_\lambda(t-s,x-y)\sigma(s,y,u^n(s,y))W(dy,ds) \\ 
&  &  \\ 
&  & +\displaystyle\int_0^t\hspace{-2mm}\displaystyle\int_{\R}G_\lambda(t-s,x-y)b(s,y,u^n(s,y))dyds,
\end{array}
\end{equation}
for all $n\geq 0$. We  start by  proving that given $t>0$, $2\le p< \infty,$

\begin{equation}
\label{borne}
\sup_{n\geq 0}\sup_{0\le s\le t}\sup_{x \in \R}E(\vert u^n(s,x)\vert^p)\le C <+\infty,
\end{equation}
where $C$ is a constant depending on $p, t,$ the supremum norm of $u_0$ and the Lipschitz constants of 
$\sigma$ and $b$. Indeed,
\begin{equation}
  \label{u0AB}
E(\left| u^{n+1}\left( t,x\right) \right| ^p) \le C
\left\{ E(\vert u^0(t,x)\vert^p) + E(\vert A_n(t,x)\vert^p)
+ E(\vert B_n(t,x)\vert^p)\right\},  
\end{equation}
where $ A_n(t,x)$ is the second term in (\ref{picard}) and $ B_n(t,x)$ 
is the third term in the right-hand side of the same equation.\\

\noindent
Then  Jensen inequality with respect to the probability  measure $G_\lambda(t, x-y) dy$ yields
$$
\vert u^0(s,x)\vert^p \le  
 \left(\int_{\R} G_\lambda(s, x-y) \vert u_0\left( y\right)\vert^p dy \right).
$$
Taking expectation and  applying Fubini's theorem  we obtain : 
$$E(\vert u^0(s,x)\vert^p) \le  \sup_{y\in \R} E(\displaystyle|u_0(y)|^p) \int_{\R} dy\, G_\lambda(s, x-y) \le  \sup_{y\in \R} E(\displaystyle|u_0(y)|^p) .
$$
Now  as ${\bf (H_1.1)}$ holds, we get :
\begin{equation}
  \label{u0}
 \sup_{0\le s\le t}\sup_{x \in \R}E(\vert u^0(s,x)\vert^p) \le C < \infty, 
\end{equation}
for some positive constant $C$.  
\medskip\noindent\\
Burkholder's inequality yields, for any $p\geq 2$ 
$$E(\left| A_n(t,x) \right|)^p\le C E\left(\displaystyle\int_0^t%
\displaystyle\hspace{-2mm}\int_{\R}G_\lambda^2(t-s,x-y)\sigma^2(s,y,u^n(s,y))\;dyds\right)^{p/2}.$$
Set
$$\nu_t=\int_0^t\hspace{-2mm}\int_{\R}G_\lambda^2(t-s,x-y)dyds,$$
 Since $\lambda > 1$, 
$\nu_t \le \int_0^T\hspace{-2mm}\int_{\R}G_\lambda^2(t-s,x-y)dyds < \infty$
 by lemma \ref{tech}(\textbf{e}).\\
Consider
\begin{equation}
\label{jts}
J(t-s)=\int_{\R}G_\lambda^2(t-s,y)dy.
\end{equation}
Due to the scaling property (see lemma \ref{tech} (\textbf{b}) ), one easily checks that 
\begin{equation}
\label{J}
J(t-s) = C (t-s)^{-1/\lambda}.
\end{equation}
 Indeed,
$J(t-s) = (t-s)^{-1/\lambda}\int_\R G_\lambda^2(1,x)dx = (t-s)^{-1/\lambda}\int_\R \exp( -\vert \xi \vert^{2\lambda}) d\xi.$
the last equality resulting from Plancherel identity.\smallskip\par\noindent
Because of the hypotheses on the coefficients  $\sigma$ and $b$, the H\"older inequality (\ref{holder}) applied with
 $f=\sigma^2(s,y,u^n(s,y))$, $h=G_\lambda^2(t-s,x-y)$ and $q=p/2$ implies
\begin{eqnarray*}
 E(\left| A_n(t,x) \right|)^p &\le& C \,\nu_t^{\frac{p}{2}-1} E\left( \displaystyle\int_0^t%
\hspace{-2mm}\displaystyle\int_{\R}G_\lambda^2(t-s,x-y)\sigma^p(s,y,u^n(s,y))\;dyds\right) \\ 
&\le& C \left( \displaystyle\int_0^t\hspace{-2mm}\displaystyle%
\int_{\R}(1+\sup_{y \in \R}E(\vert u^n(s,y)\vert^p) \,G_\lambda^2(t-s,x-y)dyds\right) \\
\le&C& \left( \displaystyle\int_0^t(1+\sup_{y \in \R}E(\vert u^n(s,y)\vert^p)
 \left(\displaystyle\int_{\R}G_\lambda^2(t-s,x-y)dy\right)ds\right).
\end{eqnarray*}
Hence
\begin{equation}
  E(\left| A_n(t,x) \right|)^p \le C \displaystyle\int_0^t \left(1+\sup_{0\le s\le t}\sup_{y \in \R}E(\vert u^n(s,y)\vert^p)\right)
 J(t-s)ds.
\label{A}
\end{equation}
The linear growth assumption on $b$ and H\"older's inequality applied to integrals with respect to the measure
$G_\lambda(t-s,x-y)dsdy$ implies
\begin{equation}
  \label{B}
 E(\left| B_n(t,x) \right|^p )\le C\, \int_0^t\,
 \left(1+\sup_{0\le s\le t}\sup_{y \in \R}E(\vert u^n(s,y)\vert^p)\right)\,ds. 
\end{equation}
 \noindent
Collecting (\ref{u0AB}),(\ref{u0}),(\ref{A}),(\ref{B}) and  (\ref{J}) we conclude that
\begin{eqnarray*}
&&E(\left| u^{n+1}\left( t,x\right) \right|^p )\\
&& \, \le C\,\left(E(\vert u^0(t,x)\vert^p) 
+ \int_0^t\left(1+\sup_{y \in \R}E(\vert u^n(s,y)\vert^p)\right)
 (J(t-s)+1)ds\right)\\
 && \, \le C\,\left(1
+ \int_0^t \,(t-s)^{-\frac{1}{\lambda}}\sup_{0\le s\le t}
\sup_{y \in \R}E(\vert u^n(s,y)\vert^p) \,ds
  \,\right).\\
 \end{eqnarray*}
Thus by   lemma \ref{gronwall} (see appendix) we obtain (\ref{borne}). \\

\noindent
In order to prove that $(u_n(t,x),\, n\geq 0)$ converges in $L^p$, let 
$n\geq 0$, $0\leq t\leq T$ and set 
$$
M_n(t) = \sup_{0\le s\le t}\sup_{x\in \R} E(\left| u^{n+1}\left( s,x\right)
-u^{n}\left( s,x\right) \right|^p). 
$$
Using the Lipschitz property of  $\sigma$ and $b$, a similar computation implies
$$M_n(t) \le C\, \int_0^t ds \,M_{n-1}(s)
 (J(t-s)+1) .$$
\noindent
Moreover,  owing to  (\ref{borne})  we have  
$$\sup_{0\le t\le T} M_0(t)\le
\sup_{0\le t\le T}\sup_{x \in \R}E(|u^1(t,x)|^p) + 
\sup_{0\le t\le T}\sup_{x \in \R}E(|u^0(t,x)|^p) <\infty.
$$ Therefore, by lemma \ref{gronwall}
the sequences  $(u_n(t,x),\, n\geq 0)$ converges in $L^p(\Omega,{\cal{G}},P)$, uniformly in $x\in \R$ and  
$0\le t \le T$,  to a limit
$u(t,x)$. It is easy to see that $u(t,x)$ satisfies (\ref{walsh}), (\ref{borne0})
which proves the existence of a solution. Following the same approach as in Walsh \cite{Wa}, we can 
prove that the process $(u(t,x),\, t\geq 0, x\in \R)$ 
has a jointly measurable version which is continuous in $L^p$ and fulfills (\ref{walsh}). 
Uniqueness of the solution 
is checked by standard arguments. \hfill $\Box$\\
\section{H\"older continuity of the solution}
\setcounter{equation}{0}
In this section we  analyze the path regularity of $u(t,x)$. The next result extends and improves
similar estimates known for the stochastic heat equation (corresponding to the case  $\lambda=2$).
\begin{Theorem}
\label{holdercont} Let $\lambda \in ]1,2]$. Suppose that ${\bf (H_0)}$,  ${\bf (H_1.1)}$ and 
${\bf (H_1.2)}$ are satisfied. 
Then,  $\omega$-almost surely, the function $\left( t,x\right)
\longmapsto u\left( t,x\right) \left( \omega \right) $ belongs to H\"older space 
${\cal{C}} ^{\alpha ,\beta }\left( \left[ 0,T\right] \times \R \right) $ 
for $0<\alpha < (\frac{\rho}{\lambda} \wedge \frac{ \lambda  -1}{2 \lambda })$ and 
$0<\beta< (\rho \wedge \frac{\lambda-1}{2})$.
\end{Theorem}
\noindent \textbf{Proof.} 
\noindent
Fix $T>0, h>0$  and $p\in ]1,1/\rho[.$  We  show first that
\begin{equation}
\label{temps}
\sup_{0\le t\le T}\sup_{x\in \R} E(\mid u(t+h,x)-u(t,x)\mid ^p)\le C\,h^{\alpha p},
\end{equation}
for any $0<\alpha < (\frac{\rho}{\lambda} \wedge \frac{ \lambda  -1}{2 \lambda })$.\\
\par\noindent
Indeed, we have
\begin{equation}
\label{sum}
E(\mid u(t+h,x)-u(t,x)\mid ^p)\le C \sum_{i=1}^4 I_i(t,h,x),
\end{equation}
where
\begin{eqnarray*}
I_1(t,h,x)&=& E\left\vert \displaystyle{\int_{\R}} (G_\lambda(t+h,x-y)-
G_\lambda(t,x-y))u_0(y)dy \right\vert^p,\\
I_2(t,h,x)&=& E\left(\left\vert\int_0^t  \displaystyle{\int_{\R}} 
[G_\lambda(t+h-s,x-y)-G_\lambda(t-s,x-y)]\right.\right.\\
&&\hspace{3cm}\times \left.\sigma(s,y,u(s,y))W(dy,ds) \Bigg\vert^p \right),\\
I_3(t,h,x)&=&E\left(\left\vert\int_t^{t+h} 
\hspace{-2mm} \displaystyle{\int_{\R}} 
G_\lambda(t+h-s,x-y)\sigma(s,y,u(s,y))W(dy,ds) \right\vert^p \right),\\
I_4(t,h,x)&=& E\left(\left\vert\int_0^{t+h} ds\displaystyle{\int_{\R}} dy\, 
G_\lambda(t+h-s,x-y)b(s,y,u(s,y))\right.\right.\\
&&\quad  \left.\left. - \int_0^t ds \displaystyle{\int_{\R}} dy\, 
G_\lambda(t-s,x-y)b(s,y,u(s,y)) \right\vert^p \right).\\
\end{eqnarray*}
Using the semigroup property of the  convolution kernel $G_\lambda$,
$$ G_\lambda(t+h ,x-y) = \int_{\R} G_\lambda(t,x-y-z)\, G_\lambda(h,z)\, dz.$$
Hence
$$
I_1(t,h,x) =E\left(\left\vert \displaystyle{\int_{\R}}G_\lambda(h,z)\left(
\displaystyle{\int_{\R}} G_\lambda(t,x-y)\left( u_0(y-z)-u_0(y)\right) dy \right)dz\right\vert^p\right).
$$
With  H\"older's inequality (\ref{holder}), the assumption 
${\bf (H_1.2)}$ and Fubini's theorem we obtain
\begin{eqnarray}
 I_1(t,h,x)& & \le  \displaystyle{\int_{\R}}\, G_\lambda(h,z) \sup_{ y\in \R} E\vert 
u_0(y-z)-u_0(y)\vert^p\, dz\nonumber\\
& & \quad \le C \displaystyle{\int_{\R}}\, G_\lambda(h,z) \, \vert z\vert^{\rho\,p} \,dz. 
\end{eqnarray}
Now, due to the self-similarity property (see lemma \ref{tech} {\bf b })
$$
\int_{\R}\, G_\lambda(h,z) \, \vert z\vert^{\rho\,p} \,dz = 
\int_{\R}\,h^{-1/\lambda} G_\lambda(1,h^{-1/\lambda}\,z)\,\vert z\vert^{\rho\,p} \,dz
$$
$$
 =  h^\frac{\rho\, p}{\lambda} \, \int_{\R} G_\lambda(1,y)\vert y \vert^{\rho\,p}\,dy. 
$$
Using the fact that $G_\lambda(1,y) \leq \frac{C}{1+y^2}$ (see lemma \ref{tech}
{\bf c}),  and that $\rho p < 1$ we obtain that 
$$
\int_{\R} G_\lambda(1,y)\vert y \vert^{\rho\,p}\,dy < \infty.
$$
Therefore we have proved that
\begin{equation}
I_1(t,h,x) \, \le \, C\, h^\frac{\rho\, p}{\lambda}.
\label{}
\end{equation}
\noindent
Bukholder  inequality, H\"older  inequality (\ref{holder}) applied to integrals with respect to the measure
$[G_\lambda(t+h-s,x-y)-G_\lambda(t-s,x-y)]^2 ds dy$, 
the growth assumption on $\sigma$ and (\ref{borne0})  yield  the following bound on $I_2$.
\begin{eqnarray*}
&& I_2(t,h,x)\le  C \left(1+\sup_{0\le s\le t}\sup_{x\in \R}E(\vert u(s,x)\vert^p )\right)\\
&& \quad \times \left(\left\vert\int_0^t  \displaystyle{\int_{\R}} 
\, [G_\lambda(t+h-s,x-y)-G_\lambda(t-s,x-y)]^2 ds dy \right\vert^{p/2}\right)\\
&&\quad \le C \left(\int_0^t\hspace{-2mm} 
\displaystyle{\int_{\R}}\left({\cal F}(e^{-(t+h-s)\mid \,\cdot\,\mid^\lambda })(y)-  
{\cal F}(e^{-(t-s)\mid \,\cdot\,\mid^\lambda })(y) \right)^2 ds dy\right)^{p/2}.
\end{eqnarray*}

Therefore, using  Plancherel identity  one easily checks that
 
\begin{eqnarray*}
&&\int_0^t\hspace{-2mm}
\displaystyle{\int_{\R}}\left({\cal F}(e^{-(t+h-s)\mid \,\cdot\,\mid^\lambda })-  
{\cal F}(e^{-(t-s)\mid \,\cdot\,\mid^\lambda }) \right)^2 (y)\, ds dy \\
&&\hspace{3cm} =\,\,
\int_0^t 
\displaystyle{\int_{\R}}\left(e^{-(t+h-s)\mid y\mid^\lambda }-  
e^{-(t-s)\mid y\mid^\lambda } \right)^2 \,ds dy\\
&& \hspace{3cm} =\,\,
\int_0^t \hspace{-2mm}
\displaystyle{\int_{\R}} \, e^{-2(t-s)\mid y\mid^\lambda}\left(e^{-h\mid y\mid^\lambda }- 1\right)^2 ds \,dy.
\end{eqnarray*}

\noindent
 Decomposing the integral on $\R$ into  integrals
 on $\{\vert y\vert>1\}$ and its complementary set, 
we have

$$I_2(t,h,x)\le C\,(I_{2,1}(t,h,x) + I_{2,2}(t,h,x))$$
\noindent
where

\begin{eqnarray*}
I_{2,1}(t,h,x)&=&\left(\int_0^t  \hspace{-2mm}
\displaystyle{\int_{\mid y\mid \le 1}} \,
e^{-2(t-s)\mid y\mid^\lambda}\left(e^{-h\mid y\mid^\lambda }- 1\right)^2 ds dy\right)^{p/2},\\
I_{2,2}(t,h,x)&=&\left(\int_0^t\hspace{-2mm}\displaystyle{\int_{\mid y\mid >1}}\,
e^{-2(t-s)\mid y\mid^\lambda}\left(
 e^{-h\mid y\mid^\lambda }- 1\right)^2 ds dy\right)^{p/2}.
\end{eqnarray*}

\noindent
Then by the mean value theorem,
\begin{eqnarray*}
\int_0^t \hspace{-2mm}
\displaystyle{\int_{\mid y\mid \le 1}}\,e^{-2(t-s)
\mid y\mid^\lambda}\left(e^{-h\mid y\mid^\lambda }- 1\right)^2 ds dy&\le&
\int_0^T \hspace{-2mm}
\displaystyle{\int_{\mid y\mid \le 1}}\,e^{-2(t-s)
\mid y\mid^\lambda} h^2 ds dy\\
&\le& C h^2.
\end{eqnarray*}

\noindent
On the set $\{\vert y\vert>1\}$, let $0<\alpha<\frac{\lambda - 1}{2 \lambda}$,
then the same argument as above implies

\begin{eqnarray*}
&&\int_0^t \hspace{-2mm}
\displaystyle{\int_{\mid y\mid >1}}e^{-2(t-s)
\mid y\mid^\lambda}\left(e^{-h\mid y\mid^\lambda }-1\right)^2dsdy\\ 
&&\hspace{1cm} =\,\,
\int_0^t \hspace{-2mm}
\displaystyle{\int_{\mid y\mid >1}} e^{-2(t-s)
\mid y\mid^\lambda}\left(1 - e^{-h\mid y\mid^\lambda}\right)^{2\alpha}
\left(1 - e^{-h\mid y\mid^\lambda }\right)^{2-2\alpha}ds dy\\
&& \hspace{1cm} \le \,\,C \int_0^\infty \hspace{-2mm}
\displaystyle{\int_{\mid y\mid >1}} \,e^{-2s
\mid y\mid^\lambda}\vert h \vert^{2\alpha} \vert  y\vert^{2\lambda\alpha} ds dy
\\
&& \hspace{1cm} \le \,\,C 
\displaystyle{\int_{\mid y\mid >1}}\,\vert h \vert^{2\alpha} \vert y \vert^{2\lambda\alpha}
\vert y\vert^{-\lambda}dy \\
&&  \hspace{1cm} \le \,\,C \,h^{2\alpha} \displaystyle{\int_{\mid y\mid >1}}\, \vert y \vert^{\lambda(2\alpha -1)}dy
\le C h^{2\alpha}.
\end{eqnarray*}

\noindent
Consequently, for $0<\alpha < \frac{\lambda - 1}{2 \lambda}$, we have proved that
$$I_{2,1}(t,h,x)\le C \, h^{ p},$$ 
$$I_{2,2}(t,h,x)\le C\, h^{\alpha p}.$$
Since $0 < \alpha  <\frac{\lambda - 1}{2 \lambda} <1,\,\forall \lambda \in ]1,2] $,
we obtain
\begin{equation}
\label{i2}
I_2(t,h,x)\le C\, h^{\alpha p}.
\end{equation}
\noindent
As before, Bukholder  inequality, H\"older  inequality (\ref{holder}) applied to integrals with respect to the measure
$G^2_\lambda(t+h-s,x-y) ds dy$, 
the growth assumption on $\sigma$ and (\ref{borne0})  yield
\begin{eqnarray*}
 I_3(t,h,x)&\le& C\left(1+\sup_{0\le s\le t}\sup_{x\in \R}E(\vert u(s,x)\vert^p )\right) \\
&&\quad\quad\times \left( \int_t^{t+h}\int_{\R}
G^2_\lambda(t+h-s,x-y)\,dsdy\right)^{p/2}.
\end{eqnarray*}
Recalling from (\ref{J}) that 
$$\displaystyle{\int_{\R}} G^2_\lambda(t+h-s,x-y)\,dy = J(t+h-s)
= C(t+h-s)^{-1/\lambda}$$
we compute $\displaystyle\int_t^{t+h} (t+h-s)^{-1/\lambda} ds  =
C\,h^{\frac{\lambda-1} {\lambda}}$.
\par\noindent Thus 
\begin{equation}
  \label{eq:I3}
  I_3(t,h,x)\le  C\,h^{\frac{p(\lambda-1)} {2\lambda}}.
\end{equation}
A change of variable yields

$$I_4(t,h,x)\le C\,(I_{4,1}(t,h,x) + I_{4,2}(t,h,x))$$

\noindent
with
\begin{eqnarray*}
I_{4,1}(t,h,x)&=& E\left(\left\vert\int_0^{h}ds\displaystyle{\int_{\R}}dy 
G_\lambda(t+h-s,x-y)b(s,y,u(s,y)) \right\vert^p \right),\\
I_{4,2}(t,h,x)&=&
E\left(\left\vert\int_0^{t} ds\displaystyle{\int_{\R}} dy\, 
G_\lambda(t-s,x-y)\right.\right.\\
&& \quad\quad\times \left.\left(b\left(s+h,y,u(s+h,y)\right)- 
b\left(s,y,u(s,y)\right)\right)\Bigg\vert^p \right).
\end{eqnarray*}
 Applying H\"older  inequality (\ref{holder})  to integrals with respect to the measure $G_\lambda(t+h-s,x-y)\,ds dy$, 
the growth assumption on $b$ and (\ref{borne0})  we get
\begin{eqnarray*}
I_{4,1}(t,h,x)&\le& C \left(1+\sup_{0\le s\le t}\sup_{x\in \R}E(\vert u(s,x)\vert^p )\right) \\
&& \quad\times \left(\displaystyle\int_0^{h} ds \displaystyle{\int_{\R}} dy\,
G_\lambda(t+h-s,x-y)\right)^{p}.
\end{eqnarray*}
Since $\int_{\R} G_\lambda(t+h-s,x-y)\,dy = 1$, we obtain
\begin{equation}
\label{i41}
I_{4,1}(t,h,x)\le C h^p.
\end{equation}
Again H\" older inequality applied to integral w.r.t. the measure
 $G_\lambda(t-s, x-y)\, dsdy$,  Fubini's theorem and the Lipschitz property of $b$ imply
\begin{eqnarray*}
I_{4,2}(t,h,x)& \le & C \left(\int_0^t \left(h^p + \sup_{y\in \R}E(\vert u(s+h,y)-u(s,y)\vert^p )\right)\, ds\right) \\
 & & \times
\left( \int_0^{T}\int_{\R}  
G_\lambda(t-s,x-y) \,dsdy\right).
\end{eqnarray*}
Hence 
\begin{equation}
  \label{i42}
 I_{4,2}(t,h,x) \le  C 
 \left( h^p +\int_0^t \sup_{y\in \R}E(\vert u(s+h,y)-u(s,y)\vert^p )\,ds\right).
\end{equation}
Then, putting together (\ref{sum})-(\ref{i42}) we obtain for $0<\alpha<\frac{\lambda - 1}{2 \lambda}$

\begin{eqnarray*}
\sup_{x\in \R}E(\vert u(t+h,x)-u(t,x)\vert^p )&\le& C\,h^{p \min(\frac{\rho}{\lambda}, \alpha)}\\
 &+&
 C \int_0^t \sup_{x\in \R}E(\vert u(s+h,x)-u(s,x)\vert^p )ds.
\end{eqnarray*}
Finally, the estimates (\ref{temps}) follows from standard Gronwall's Lemma.
\newpage
\noindent
Consider now the increments in the space variable. We want to check that for any $T>0$, $p \in [2, \infty), 
x \in \R$, $z$ in a compact set $K$ of $\R$ and $\beta \in(0, \rho \wedge (\frac{\lambda-1}{2}))$,

\begin{equation}
\label{espace}
\sup_{0\le t\le T}\sup_{x\in \R} E(\mid u(t,x+z)-u(t,x)\mid ^p)\le C\,z^{\beta p},
\end{equation}

\noindent
We write

\begin{equation}
\label{sum2}
E(\mid u(t,x+z)-u(t,x)\mid ^p)\le C \sum_{i=1}^3 J_i(t,z,x),
\end{equation}
with
\begin{eqnarray*}
J_1(t,z,x)&=& E \left\vert \displaystyle{\int_{\R}} (G_\lambda(t,x+z-y)-G_\lambda(t,x-y))u_0(y)dy \right\vert ^p,\\
J_2(t,z,x)&=& E\left(\left\vert\int_0^t \hspace{-2mm} \displaystyle{\int_{\R}} 
[G_\lambda(t-s,x+z-y)-G_\lambda(t-s,x-y)]\right.\right.\\
&& \quad \quad \quad \quad \hspace{2cm}
\left.\times \sigma(s,y,u(s,y))W(dy,ds) \Big \vert^p \right),\\
J_3(t,z,x)&=& E\left(\left\vert\int_0^{t} ds \displaystyle{\int_{\R}} dy\,
[G_\lambda(t-s,x+z-y)-G_\lambda(t-s,x-y)]\right.\right.\\
&& \quad \quad \quad \quad\hspace{3,5cm}\left.\times
\,b(s,y,u(s,y))\Big\vert^p \right).
\end{eqnarray*}
In the remainder of the proof we are going to establish separate upper bounds for $J_1, J_2$ and $J_3$.\\

\noindent
A change of variable gives immediately
$$
J_1(t,z,x) =  E \left\vert \displaystyle{\int_{\R}} G_\lambda(t,x-y)\left(u_0(y+z)
-u_0(y)\right)\,dy \right\vert ^p.
$$
Applying again  H\"older's inequality (\ref{holder})  to integral w.r.t. the measure
$G_\lambda(t,x-y)\,dy$, the assumption 
${\bf (H_1.2)}$ and Fubini's theorem we obtain
\begin{eqnarray*}
J_1(t,z,x)&\le & \, C \left(\displaystyle{\int_{\R}}\, G_\lambda(t,x-y) \sup_{ y\in \R} E(\vert 
u_0(y+z)-u_0(y)\vert^p)\, dy\right)\\
&\le & \, C \left(\displaystyle{\int_{\R}}\, G_\lambda(t,x-y) \vert z\vert^{\rho\,p}\, dy\right)\le \, C 
\,\vert z\vert^{\rho\,p}.
\end{eqnarray*}

\noindent
Bukholder's inequality and H\"older's inequality (\ref{holder}) applied to integrals with respect to the
measure $[G_\lambda(t-s,x+z-y)-G_\lambda(t-s,x-y)]^2 ds dy$, the linear growth assumption on $\sigma$ 
and (\ref{borne0}) imply
\begin{eqnarray*}
J_2(t,z,x)&\le & C \,\left(1+\sup_{0\le t\le T}\sup_{x\in \R}E(\vert u(t,x)\vert^p )\right)\\
&&  \times \left(\left\vert\int_0^t ds \displaystyle{\int_{\R}} 
dy\; [G_\lambda(t-s,x+z-y)-G_\lambda(t-s,x-y)]^2 \right\vert^{p/2} \right)\\
J_2(t,z,x)&\le & C\left(\int_0^t  
\int_{\R}  \left({\cal F}(e^{-2\pi i z \cdot}e^{-(t-s)\mid\cdot\mid^\lambda })
-  {\cal F}(e^{-(t-s)\mid \cdot\mid^\lambda }) \right)^2 (x-y)\right)^{p/2}\\
&\le & C \,\left(\int_0^t ds
\displaystyle{\int_{\R}}dy \left(e^{-2\pi i z y}e^{-(t-s)\mid \,y\,\mid^\lambda }-  
e^{-(t-s)\mid \,y\,\mid^\lambda } \right)^2 \right)^{p/2}\\  
&\le & C\,( J_{2,1}(t,z,x)+ J_{2,2}(t,z,x)), 
\end{eqnarray*}
\noindent
where we have used the property that 
${\cal F}\left(f(x)\right)(\xi + a) = 
{\cal F}\left( e^{-2i\pi ax} f(x)\right)(\xi)$
 and the Plancherel identity and denote 
\begin{eqnarray*}
J_{2,1}(t,z,x)& = & \left(\int_0^t ds
\displaystyle{\int_{\vert y \vert \le 1}}dy \left(e^{-2\pi i z y}e^{-(t-s)\mid \,y\,\mid^\lambda }-  
e^{-(t-s)\mid \,y\,\mid^\lambda } \right)^2 \right)^{p/2},\\
J_{2,2}(t,z,x)& = & \left(\int_0^t ds
\displaystyle{\int_{\vert y \vert >1}}dy \left(e^{-2\pi i z y}e^{-(t-s)\mid \,y\,\mid^\lambda }-  
e^{-(t-s)\mid \,y\,\mid^\lambda } \right)^2 \right)^{p/2}.
\end{eqnarray*}

\noindent
We therefore have, by the mean value theorem

\begin{equation}
\label{j21}
J_{2,1}(t,z,x) \le C \vert z\vert^p.
\end{equation}

\noindent
On the other hand,  for any $0<\beta <\frac{\lambda - 1}{2}$

\begin{eqnarray}
\label{j22}
&&J_{2,2}(t,z,x)\nonumber\\
&& \quad =\left(\int_0^t\hspace{-2mm}
\displaystyle{\int_{\vert y \vert >1}} e^{-(t-s)\mid y\mid^\lambda }
\left(e^{-2\pi i z y}-  
1 \right)^{2 \beta} \left(e^{-2\pi i z y}-  
1 \right)^{2 - 2 \beta}ds dy\right)^{p/2}\nonumber\\
& & \quad \le C \left(\int_0^t ds
\displaystyle{\int_{\vert y \vert >1}}dy\,  e^{-(t-s)\mid \,y\,\mid^\lambda }
\vert y\vert^{2 \beta} \vert z\vert^{2 \beta}\right)^{p/2}\nonumber\\
& & \quad \le C \vert z\vert^{2 \beta} \left(\displaystyle{\int_{\vert y \vert >1}}dy\, \vert y\vert^{2 \beta}  
\int_0^t ds\, \, e^{-(t-s)\mid \,y\,\mid^\lambda }
 \right)^{p/2} \nonumber\\
&& \quad \le C \vert z\vert^{ \beta p} \displaystyle{\int_{\vert y \vert >1}} 
\frac{dy}{\vert y\vert^{\lambda - 2 \beta}} \le C \vert z\vert^{\beta p}.
\end{eqnarray}
Finally, by a change of variable, the Lipschitz property of $b$, and H\"older's inequality,
\begin{eqnarray}
\label{j3}
J_3(t,z,x)&\le & E\left(\left\vert\int_0^{t} ds \displaystyle{\int_{\R}} dy\,
G_\lambda(t-s,x-y)\right.\right.\nonumber\\
&& \quad\quad \times\Bigg.\left.
[b(s,y+z,u(s,y+z))- b(s,y,u(s,y))]\right\vert^p \Bigg)\nonumber\\
&\le &  
C\left( z^p + \int_0^t \sup_{y\in \R}E(\vert u(s,y+z)-u(s,y)\vert^p )\,ds\right). 
\end{eqnarray}

\noindent 
Then (\ref{espace}) follows from (\ref{sum2})-(\ref{j3}) and Gronwall's lemma.
The H\"older continuity in the time and space variables results from Kolmogorov criterion.\hfill $\Box$

\begin{remark}
Hypothesis  ${\bf (H_1.2)}$ is useful to have H\"older continuity up to time 0.
 If we discard ${\bf (H_1.2)}$ and assume instead that 
$$
 {\bf (H_1.3)}\quad E\left(\int_{\R}\displaystyle|u_0(y)|\,dy\right)^p < \infty
$$ then  $\omega$-almost surely, the function $\left( t,x\right)
\longmapsto u\left( t,x\right) \left( \omega \right) $ belongs
to
${\cal{C}} ^{\alpha ,\beta }\left( \left[ \epsilon,T\right] \times \R \right) $ for $0<\alpha <  \frac{ \lambda  -1}{2 \lambda }$ and 
$0<\beta< \frac{\lambda-1}{2}$, for any $\epsilon >0$.
\end{remark}

Indeed, we slightly modify the preceding proof to bound 
$$I_1(t,h,x)= E\left\vert \displaystyle{\int_{\R}} (G_\lambda(t+h,x-y)-
G_\lambda(t,x-y))u_0(y)dy \right\vert^p,
$$
and
$$
J_1(t,z,x) =  E \left\vert \displaystyle{\int_{\R}} (G_\lambda(t,x+z-y)-G_\lambda(t,x-y))u_0(y)dy \right\vert ^p.
$$ 

First we bound
$$I_1(t,h,x) \leq
 \sup_{z\in \R}\vert G_\lambda(t+h, z)-G_\lambda(t, z) \vert^p\,\cdot
  E\left(\int_{\R}\displaystyle|u_0(y)|\,dy\right)^p.
$$
The following estimates are elementary
$$G_\lambda(t+h, z)-G_\lambda(t, z) = 
 \int_{\R} e^{2 i\pi  z\xi} 
(e^{-(t+h)\vert\xi\vert^\lambda} - e^{-t\vert\xi\vert^\lambda})\, d\xi,
$$
$$|\,G_\lambda(t+h, z)-G_\lambda(t, z)| \leq 
\int_{\R} e^{-t\vert\xi\vert^\lambda} |\,e^{-h\vert\xi\vert^\lambda}-1|\, d\xi, 
$$
$$
 |\,G_\lambda(t+h, z)-G_\lambda(t, z)| \leq 
h \,\int_{\R} e^{-\epsilon\vert\xi\vert^\lambda}\vert\xi\vert^\lambda\, d\xi = C\,h. 
$$
Hence
 $$I_1(t,h,x) \leq C\, h^p.$$

As for the space increments, we bound

$$G_\lambda(t, x+z)-G_\lambda(t, x) = 
 \int_{\R} e^{2 i\pi  x \xi} e^{-t\vert\xi\vert^\lambda}
(e^{2 i\pi  z \xi} - 1)\, d\xi,
$$

$$
 |\,G_\lambda(t, x+z)-G_\lambda(t, x)| \leq 
 \vert z\vert\,
\int_{\R} e^{-\epsilon\vert\xi\vert^\lambda}2 \pi \vert\xi\vert\, d\xi =
 C\,\vert z\vert, 
$$
Hence
 $$J_1(t,z,x) \leq C\, z^p.$$
The rest of the proof is the same as for theorem \ref{holdercont}. \hfill $\Box$
\section{Appendix}
\setcounter{equation}{0}

\begin{Lemma}
Let $f, h$ be two functions defined on $\R$ and $\mu$ a positive measure such that $f \cdot h \in L^1(\mu)$. Then, for all
$q>1$, we have:
\begin{equation}
\label{holder}
\left\vert\int  f\cdot \vert h \vert d\mu \right\vert^q\le \left(\int \vert f\vert^q\cdot \vert h\vert d\mu\right) 
\left(\int \vert h\vert d\mu\right)^{q-1}.
\end{equation}
\end{Lemma}
\noindent \textbf{Proof.} Set $\nu= \vert h\vert d\mu$, then the result follows from the H\" older inequality applied to
$\int fd\nu$. \hfill $\Box$
\smallskip\par\noindent
The following elementary Lemma is an extension of Gronwall's Lemma  akin to   lemma 3.3  established in Walsh \cite{Wa}.
\begin{Lemma} 
\label{gronwall}Let $\theta > 0$.
 Let $(f_n, n\in \N)$
be a sequence of non-negative functions on $[0,T]$ and $\alpha, \beta$ be  non-negative real  numbers such that for $0\le t\le T,\, n\geq 1$
\begin{equation}
f_n(t)\le \alpha+ \int_0^t \beta\,f_{n-1}(s) (t-s)^{\theta -1}ds.
\label{gronvol}
\end{equation}
If $\sup_{0\le t \le T}f_0(t)=M, $ then for $n\geq 1$,
\begin{equation}
f_n(t)\le \frac{1}{2} \left(\alpha + \alpha
\exp  \left(\frac{2 \beta t^{\theta}}{\theta}\right)  +
 \frac{M}{n!} \left( \frac{2 \beta t^\theta}{\theta}\right)^n\right).
\end{equation}
In particular, $\sup_{n \geq 0}\sup_{0\le t \le T}f_n(t)<\infty $, and if $\alpha=0$, then
$\sum_{n \geq 0} f_n(t)$ converges uniformly on $[0,T]$.
\end{Lemma}
\par\noindent \textbf{Proof.} 
Let us prove  by induction that, for $n\geq 1$,
\begin{equation}
\label{induc}
f_n(t)\le \alpha \Bigg(1 + \sum_{1\leq k \leq n-1} \frac{2^{k-1}}{k!}
\left(\frac{\beta t^\theta}{\theta}\right)^k \Bigg) + \;M\,
 \frac{2^{n-1}}{n!} \left( \frac{\beta t^\theta}{\theta}\right)^n.
\end{equation}
 The initial step is readily checked :
$$
f_1(t) \leq \alpha + \int_0^t \beta\,M (t-s)^{\theta -1}ds
=\alpha + M\,\frac{\beta t^\theta}{\theta}.
$$
 Now since \ref{gronvol}
we have 
$$
f_n(t)\le \alpha +
\int_0^t \beta \left( \alpha + \alpha
  \sum_{1\leq k \leq n-2} \frac{2^{k-1}}{k!}
\left(\frac{\beta s^\theta}{\theta}\right)^k \right.
$$
\begin{equation}
+  \left.  M\,\frac{ 2^{n-2}}{(n-1)!}
 \left( \frac{\beta s^\theta}{\theta}\right)^{n-1}\right) 
(t-s)^{\theta -1}\,ds.
\end{equation}
Consider
\begin{equation}
   \int_0^t s^{k\theta} (t-s)^{\theta-1} ds \leq
  \int_0^{t/2}  (t-s)^{(k+1)\theta-1} ds + 
\int_{t/2}^t s^{(k+1)\theta -1} ds.
\end{equation}
%
Hence we may bound
\begin{equation}
  \label{eq:beta}
   \int_0^t s^{k\theta} (t-s)^{\theta-1} ds 
 \leq 2 \frac{t^{(k+1)\theta}}{(k+1)\theta}. 
\end{equation}
Summation over $k$ brings (\ref{induc}). \hfill $\Box$ 
{}


\addcontentsline{toc}{section}{References}
\begin{thebibliography}{99}
\bibitem{Ba}  Bardos, C.; Penel, P.; Frisch, U.; Sulem P.L. (1979) Modified dissipativity for a nonlinear
evolution equation arising in turbulence. \textit{Arch. Rat. mech. anal.} \textbf{71},
237--256.

\bibitem{Bi}  Biller, P. and Woyczynski, W.A. (1998) Global and exploding solutions for nonlocal quadratic 
evolution problems. \textit{SIAM J. Appl. Math.} \textbf{59}%
, 845-869.

\bibitem{Ca}  Cairoli, R. and Walsh, J. B. (1975) Stochastic integrals in the plane. \textit{Acta Math.} 
\textbf{134}, 111--183.


\bibitem{Da}  Dawson, D.A. and Gorostiza L. G. (1990) Generalized solutions of a class of nuclear-space-valued
stochadtis evolutions equations. \textit{Appl. Math. Optim.} \textbf{22},
241--263.

\bibitem{Do}  Donati--Marin, C. and Pardoux, E. (1993) White noise driven
SPDEs with reflection. \textit{Probab. Theory Relat. Fields}. \textbf{95},
1--24.

\bibitem{Dro1}  Droniou, J.; Gallouet, T; Vovelle, J. (2003)  Global solution and smoothing 
effect for a non-local regularization of a hyperbolic equation. 
\textit{J. Evol. Equ.} \textbf{3}, no. 3, 499--521.

\bibitem{Dro}  Droniou, J. and Imbert, C. (2004)  Fractal first order partial differential equations. 
\textit{Preprint}.

\bibitem{It}  It\^o, K. (1984) Lectures on stochastic Processes. Springer, Berlin. 



\bibitem{Ko}  Komatsu, T. (1984) On the martingale problem for genarators of stable processes with perturbations.
 \textit{Osaka J. Math}. \textbf{21}, 113--132.


\bibitem{Le}  L\'evy, P. (1925) Calcul des probabilit\'es.

\bibitem{Man}  Mann, J.A. and  Woyczynski, W.A. (2001) Growing fractal interfaces in the presence
 of self-similar hopping
surface diffusion. \textit{Physica A}. \textbf{291}, 159-183.


\bibitem{Me}  Metzler, R., Klafter J. (2000) The random walk's guide to anomalous diffusion: a 
fractional dynamics approach. \textit{Physics reports}. \textbf{339}, 281--292.
1-77.



\bibitem{Pa}  Pardoux, E. (1979) Stochastic partial differential equations and
filtering of diffusion processes. {Stochastics}. \textbf{3}, 127--167.

\bibitem{Pa2}  Pardoux, E. (1993) Stochastic partial differential equations, a
review. \textit{Bull. Scien. Math}. II s\'{e}rie \textbf{117}, 29--47.

\bibitem{Po}  Podlubny, I. (1999) \textit{Fractional Differential equations:
an Introduction to Fractional Derivatives, Fractional Differential equations, to Methods of Their Solution and Some of their Applications.}%
Academic Press, San Diego, CA.


\bibitem{St}  Stroock, D.W. (1975), Diffusion processes associated with L\'evy generators. \textit{Z. Wahr. verw.
Geb.}. \textbf{32}, 209-244.

\bibitem{Wa}  Walsh, J.B. (1986) An Introduction to stochastic partial
differential equations. Ecole d'Et\'{e} de Probabilit\'{e}s de Saint--Flour
XIV 1984. \textit{Lecture Notes in Mathematics.} \textbf{1180}, pp.
265--437, Berlin, Heidelberg, New York, Tokyo, Springer--Verlag.

\end{thebibliography}
\end{document}